\documentclass[preprint,12pt]{elsarticle}
\usepackage{amsthm,amsmath,amssymb}
\usepackage[colorlinks=true,citecolor=black,linkcolor=black,urlcolor=blue]{hyperref}
\usepackage{graphicx}
\usepackage{float}
\usepackage{epstopdf}
\usepackage{epsfig}
\usepackage{subfigure}

\theoremstyle{plain}
\newtheorem{theorem}{Theorem}
\newtheorem{lemma}[theorem]{Lemma}

\theoremstyle{definition}
\newtheorem{definition}[theorem]{Definition}
\newtheorem{example}[theorem]{Example}

\theoremstyle{remark}
\newtheorem{remark}[theorem]{Remark}

\title{Checkerboard colourable twisted duals}
\author{Qi Yan~~~Xian'an Jin\footnote{Corresponding author.}\\
\small School of Mathematical Sciences\\[-0.8ex]
\small Xiamen University\\[-0.8ex]
\small P. R. China\\
\small\tt Email:qiyanmath@163.com;xajin@xmu.edu.cn}
\date{}

\begin{document}

\begin{abstract}
In this note we show that any embedded graph has a checkerboard colourable twisted dual and any Eulerian embedded graph has a checkerboard colourable partial Petrial, answering questions posed by Ellis-Monaghan and Moffatt. The proofs are based on orientations of their medial graphs and orientations of boundary components of their corresponding ribbon graphs. The arrow presentations of ribbon graphs are also used. We also obtain two related results.
\end{abstract}
\maketitle

\section{Introduction}

A \emph{cellularly embedded graph} is a graph $G$ embedded in a closed surface $\Sigma$ such that every connected component of $\Sigma-G$ is a 2-cell, called a \emph{face} of the cellularly embedded graph.
We use the term embedded graph loosely to mean any of the following three equivalent representations of graphs in surfaces: cellularly embedded graphs, ribbon graphs and arrow presentations. We refer the reader to \cite{EM,M} for details and shall move from one to another freely.

A \emph{checkerboard colouring} of an embedded graph is an assignment of the colour red or colour blue to each face such that adjacent faces receive different colours. An embedded graph is said to be \emph{checkerboard colorable} if it has a checkerboard coloring. An embedded graph is said to be \emph{even-face} if the degree of each of its faces is even. A graph is said to be \emph{Eulerian} if the degree of each of its vertices is even. A graph is said to be \emph{bipartite} if it does not contain cycles of odd lengths. The following relations shown in Figure \ref{1} are well known.

\begin{figure}[htbp]\label{1}
\begin{center}
\includegraphics[width=8cm]{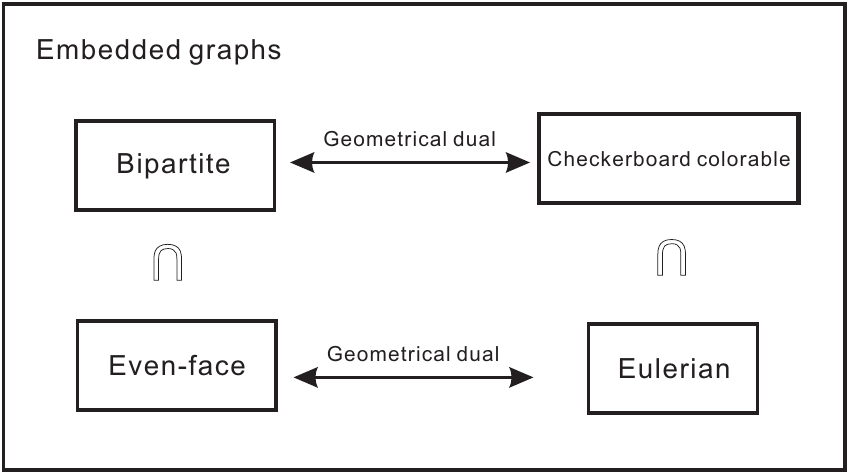}
\caption{Relations among bipartite, even-face, checkerboard colorable and Eulerian embedded graphs.}
\label{1}
\end{center}
\end{figure}

Let $G=(V(G),E(G))$ be an embedded graph. We denote by $G^{*}$ and $G^{\times}$ the \emph{geometric dual} and \emph{Petrial} \cite{WI} of $G$, respectively. Let $A\subseteq E(G)$. We denote by $G^{\delta(A)}$ (abbr. $G^{A}$) and $G^{\tau(A)}$ the partial dual \cite{chum} and partial Petrial of $G$ with respect to $A$, respectively. Partial duality and partial Petriality are further combined together to form twisted duality \cite{EM1}. We refer the reader to the monograph \cite{EM} for the details.
In \cite{EM1}, Ellis-Monaghan and Moffatt asked the following question.
\vskip0.2cm

\noindent {\bf Question 1.}
Is it possible to characterize those embedded graphs, without degree restrictions, that have a checkerboard colourable twisted dual?
\vskip0.2cm

Clearly, every bipartite embedded graph does (by taking its geometric dual). Furthermore, if $G$ and $G^{\times}$ are both orientable, then $G$ must be bipartite and hence $G^{\ast}$ is checkerboard colourable. In this note, we first prove the following theorem.

\begin{theorem}\label{th6}
Any embedded graph has a checkerboard colourable twisted dual.
\end{theorem}

It is a well-known fact that for a plane graph $G$, $G$ is checkerboard colourable if and only if $G$ is Eulerian. This fact only holds for general embedded graphs in one direction: any checkerboard colourable embedded graph is Eulerian.
And it is not true in other direction in general. In \cite{EM1}, Ellis-Monaghan and Moffatt also asked:

\vskip0.2cm
\noindent {\bf Question 2.}
If $G$ is a 4-regular embedded graph, which of its twisted dual are also 4-regular and checkerboard colorable?
\vskip0.2cm

In this note we proved:
\begin{theorem}\label{th4}
Any Eulerian embedded graph has a checkerboard colourable partial Petrial.
\end{theorem}

We also obtain the following two related results.

\begin{theorem}\label{th3}
Let $G$ be an embedded graph. If $G^{\times}$ is orientable, then $G^{\ast}$ is an Eulerian graph.
\end{theorem}

\begin{theorem}\label{th1}
Let $G$ be an embedded graph and $A\subseteq E(G)$. If $G^{A}$ is checkerboard colourable, then $G-{A}$ and $G^{\ast}-{A^{c}}$ are both checkerboard colourable (hence, also Eulerian).
\end{theorem}

\section{Preliminaries}

In this section, we provide some necessary preliminaries.
A ribbon graph results naturally from a classical cellularly embedded graph by taking a small neighbourhood and is formally defined as follows.

\begin{definition}[\cite{BR}]
A {\it ribbon graph} $G$ is a (orientable or non-orientable) surface with boundary,
represented as the union of two sets of topological discs, a set $V(G)$ of vertices, and a set $E(G)$ of edges,
subject to the following restrictions.
\begin{enumerate}
\item[(1)] the vertices and edges intersect in disjoint line segments, we call them {\it common line segments} as in \cite{Metrose};
\item[(2)] each such common line segment lies on the boundary of precisely one vertex and precisely one edge;
\item[(3)] every edge contains exactly two such common line segments.
\end{enumerate}
\end{definition}

\begin{figure}[!htbp]
\begin{center}
\includegraphics[width=15cm]{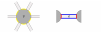}
\caption{Common line segments (red), edge line segments
(blue) and vertex line segments (yellow).}
\label{2}
\end{center}
\end{figure}

\begin{figure}[!htbp]
\begin{center}
\includegraphics[width=15cm]{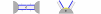}
\caption{Four half-edge line segments $\{1,2,3,4\}$of $e$, one of vertex line segments of $v$ and its incident two half-edge line segments.}
\label{3}
\end{center}
\end{figure}

We denote by $d{(v)}$ the degree of the vertex $v$ in $G$, i.e. the number of half-edges incident with $v$. Let $G$ be a ribbon graph, $v\in V(G)$ and $e\in E(G)$. By deleting the common line segments from the boundary of $v$,
we obtain $d(v)$ disjoint line segments, called \emph{vertex line segments} \cite{Metrose}.
By deleting common line segments from the boundary of $e$, we obtain two disjoint line segments,
called \emph{edge line segments} \cite{Metrose}. See Figure \ref{2}.
We assume that each edge line segment consists of two \emph{half-edge line segments}.
It is obvious that every edge disc contains four half-edge line segments.
For any vertex line segment, there are exactly two half-edge line segments incident with it as shown in Figure \ref{3}.

Let $G$ be a ribbon graph and $A\subseteq E(G)$. Then the \emph{partial Petrial},
$G^{\tau(A)}$, of $G$ with respect to $A$ is the ribbon graph obtained from $G$ by adding a half-twist to each of the edges in $A$. In \cite{chum}, Chumtov introduced the arrow presentation of a ribbon graph and defined the partial dual in terms of it.

\begin{definition}[\cite{chum}]\label{de1}
An \emph{arrow presentation} consists of a set of circles, each with a collection of disjoint, labelled arrows on them. Each label appears on
precisely two arrows.
\end{definition}

\begin{figure}[!htbp]
\centering
\includegraphics[width=14cm]{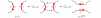}
\caption{Taking the partial dual of an edge $e$ in an arrow presentation.}
\label{4}
\end{figure}

The arrow presentation of the \emph{partial dual} $G^{A}=G^{\delta(A)}$ of $G$ with respect to $A$ is
obtained as follows. For each $e\in A$, suppose $\alpha$ and $\beta$ are the two arrows labelled $e$ in the arrow presentation of $G$.
Draw a line segment with an arrow on it directed from the head of $\alpha$ to the tail of $\beta$, and a line segment with an arrow on it directed
from the head of $\beta$ to the tail of $\alpha$. Label both of these arrows $e$ and delete $\alpha$ and $\beta$ and the arcs underlining them. See Figure \ref{4}. $G^A$ then can be recovered from its arrow presentation.

Let $\mathcal{B}=<\delta, \tau| \delta^{2}, \tau^{2}, (\delta\tau)^{3}>$.

\begin{definition}[\cite{EM1}]
Let $G$ be a ribbon graph. The ribbon graph $H$ is called the twisted dual of $G$ if it can be written in the form
$$H=G^{\prod_{i=1}^6 \xi_i(A_i)},$$
where the $A_i$'s partition $E(G)$ and the $\xi_i$'s are the six elements of $\mathcal{B}$.
\end{definition}

\begin{figure}[!htbp]
\begin{center}
\includegraphics[width=12cm]{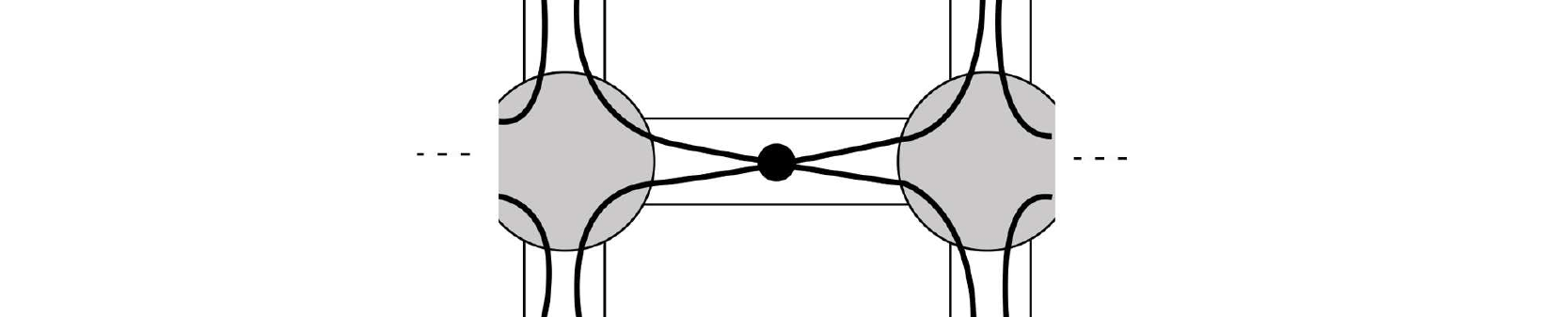}
\caption{Medial graph of a ribbon graph.}
\label{5}
\end{center}
\end{figure}

Medial graphs as a tool will be used very often throughout this note. We can form the medial graph of a ribbon graph inside the ribbon graph as shown in Figure \ref{5}. In particular, the medial
graph of an isolated vertex is a circle inside and along the boundary of the vertex disk, called \emph{free loop}.
An all-crossing direction of the medial graph $G_{m}$ of a ribbon graph $G$ is an assignment
of a direction to each edge of $G_{m}$ in such a way that at each vertex $v$ of $G_{m}$,
when we follow the cyclic order of the directed edges incident to $v$, we find
head, head, tail and tail. If $G_{m}$ is equipped with an all-crossing direction, then we can partition the vertices of $G_{m}$ into {\it $c$-vertices} and {\it $d$-vertices} according to the direction as shown in Figure \ref{9}. Accordingly edges of $G$ are divided into {\it $c$-edges} and {\it $d$-edges}.
All-crossing directions are generalized to semi-crossing orientations to solve the characterization of Eulerian partial duals of plane graphs in \cite{MJ}.

\section{Proofs}

\begin{figure}[!htbp]
\begin{center}
\includegraphics[width=15cm]{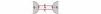}
\caption{The signs of edge-line segments and common line segments.}
\label{6}
\end{center}
\end{figure}

Given an orientable ribbon graph $G$, we fix an orientation $\mathcal{GO}$ of $G$. Let $A\subseteq E(G)$. Consider an orientation $\mathcal{BO}$ for every boundary component of $G-A$, we say that an edge-line segment of edges in $A^{c}$ or a common line segment of edges in $A$ is positive if $\mathcal{BO}$ is consistent with $\mathcal{GO}$ on this line segment and negative otherwise. See Figure \ref{6}.
Since $G-A$ is orientable, line segments on a boundary component of $G-A$ must be either all positive or all negative. Accordingly, we call the boundary component positive or negative, respectively.

\begin{figure}[!htbp]
\begin{center}
\includegraphics[width=15cm]{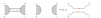}
\caption{Ribbon graphs $G, G-e$ and the arrow presentation of $G^{e}$.}
\label{7}
\end{center}
\end{figure}

\begin{lemma}\label{th5}
Let $G$ be an oriented ribbon graph and $A\subseteq E(G)$. Then there exists an orientation for every boundary component of $G-A$
such that one of two edge-line segments of $e$ or two common line segments of $f$ is positive and the other is negative for any $e\in A^{c}$ and $f\in A$
if and only if $G^{A}$ is checkerboard colourable.
\end{lemma}
\noindent {\bf Proof.} Note that $G^A$ is also orientable and there is a correspondence between boundary components of $G-A$ and those of $G^{A}$ as shown in Figure \ref{7}. We assign a face of $G^A$ blue (resp. red) if its corresponding boundary component of $G-A$ is negative (resp. positive).
The condition one of two edge-line segments of $e$ or two common line segments of $f$ is positive and the other is negative for any $e\in A^{c}$ and $f\in A$ means that adjacent faces of $G^A$ will receive different colours, completing the proof.

\begin{figure}[!htbp]
\centering
\subfigure[]{\includegraphics[width=6in]{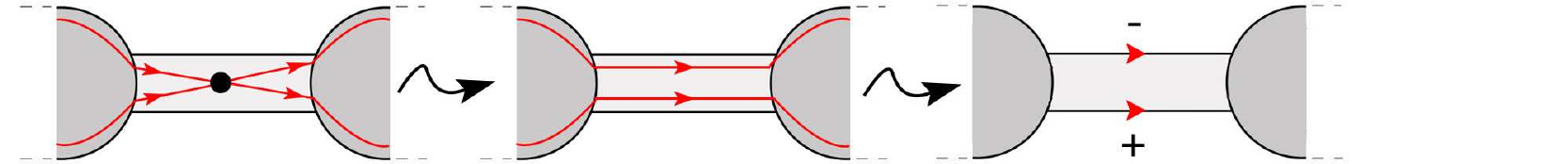}}\\
\subfigure[]{\includegraphics[width=6in]{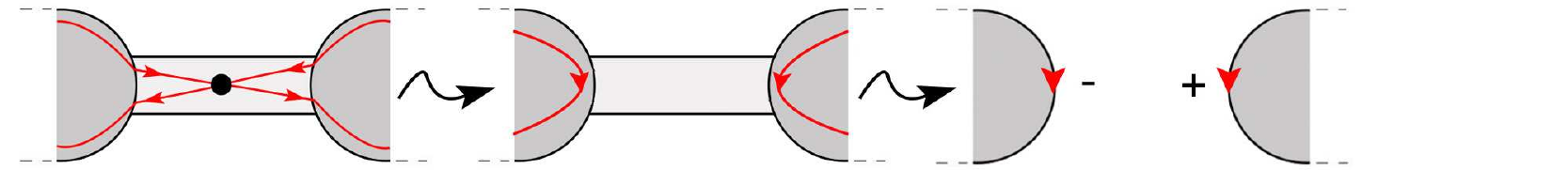}}
\caption{(a) $c$-vertex, $C$-smoothing and the signs and directions of the edge-line segments, (b) $d$-vertex, $D$-smoothing and the signs and directions of the common line segments.}\label{8}
\end{figure}

\noindent {\bf Proof of Theorem \ref{th6}.}
For any ribbon graph $G$, there exists a subset $A\subseteq E(G)$ such that $G^{\tau (A)}$ is orientable. Let $D$ be the set of $d$-edge of $G^{\tau (A)}$ arising from an all-crossing direction of $(G^{\tau (A)})_{m}$.
For any $c$-vertex (resp. $d$-vertex) of $(G^{\tau (A)})_{m}$, we do a $C$-smoothing (resp. $D$-smoothing).
This will induce an orientation of boundary components of $G^{\tau (A)}-D$ satisfying the conditions of Lemma \ref{th5}. See Figure \ref{8}. Thus $(G^{\tau (A)})^{\delta(D)}$ is checkerboard colourable.

According to the proof, for any ribbon graph $G$, we can obtain a checkerboard colourable twisted dual of $G$ as follows.

\noindent{\bf Step 1.} Draw $G$ in the plane such that each vertex disk of $G$ is not twisted and each edge disk of $G$ is twisted at most once.

\noindent{\bf Step 2.} Let $A$ be the set of all twisted edges of $G$. Convert each edge in $A$ to be non-twisted one, we obtain the orientable ribbon graph $G^{\tau(A)}$.

\noindent{\bf Step 3.} Draw the medial graph $(G^{\tau(A)})_m$ of $G^{\tau(A)}$ inside $G^{\tau(A)}$.

\noindent{\bf Step 4.} Give an all-crossing orientation of $(G^{\tau(A)})_m$ by going along straight-ahead walks of $(G^{\tau(A)})_m$.

\noindent{\bf Step 5.} Take $D$ to be the set of $d$-edges of $G^{\tau(A)}$ (i.e. $d$-vertices of $(G^{\tau(A)})_m$).

Then $(G^{\tau(A)})^{\delta(D)}$ is a checkerboard colourable twisted dual of $G$.

\vskip0.2cm

\begin{figure}[!htbp]
\begin{center}
\includegraphics[width=10cm]{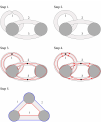}
\caption{Step 1. A ribbon graph $G$ with only one boundary component (hence, it is not checkerboard colorable). Step 2. $A=\{1\}$ and $G^{\tau(A)}$. Step 3. Medial graph $(G^{\tau(A)})_m$.
Step 4. An all-crossing orientation of $(G^{\tau(A)})_m$ and $D=\{1\}$. Step 5. $(G^{\tau(A)})^{\delta(D)}$ with two boundary components (red and blue, hence checkerboard colorable).}
\label{9}
\end{center}
\end{figure}

\begin{example}
An example is given in Figure \ref{9}.
\end{example}

\vskip0.2cm

\begin{figure}[!htbp]
\begin{center}
\includegraphics[width=13cm]{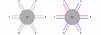}
\caption{Checkerboard colouring of half-edge vertex line segments at a vertex.}
\label{10}
\end{center}
\end{figure}

\begin{figure}[!htbp]
\begin{center}
\includegraphics[width=13cm]{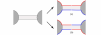}
\caption{(a) consistent edge (b) inconsistent edge.}
\label{11}
\end{center}
\end{figure}

\noindent {\bf Proof of Theorem \ref{th4}.}
Let $G$ be an Eulerian ribbon graph and $v\in V(G)$. We can use red and blue to color all half-edge line segments adjacent to $v$ and
vertex line segments of $v$ such that the colours of vertex line segments are alternating between red and blue
in the vertex boundary of $v$ and for each vertex line segment and its adjacent two half-edge line segments use the same colour as the vertex line segement.
We shall call this colouring a checkerboard colouring of $v$. See Figure \ref{10} for an example.
Now we give each vertex $v$ of $G$ a checkerboard colouring, then the edges of $G$ will be divided into consistent edges and inconsistent edges as shown in Figure \ref{11}.
Let $I$ be the set of all inconsistent edges of $G$. We consider the partial Petrial $G^{\tau(I)}$ of $G$. We could give each vertex of $G^{\tau(I)}$ the same checkerboard colouring as given in $G$. Then all edges of $G^{\tau(I)}$ are consistent. Assign the face of $G^{\tau(I)}$ the same color as its boundary, we obtain a checkerboard colouring of $G^{\tau(I)}$. Thus $G^{\tau(I)}$ is a checkerboard colourable partial Petrial of $G$.

\begin{figure}[!htbp]
\begin{center}
\includegraphics[width=15cm]{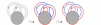}
\caption{Convert a non-checkerboard colorable Eulerian embedded graph to a checkerboard colorable one by adding half twists.}
\label{12}
\end{center}
\end{figure}

\begin{example}
The graph $G$ consisting of a vetex and two loops (a latitude circle and a longitude circle) in torus is not checkerboard colorable and it can be converted to be a checkerboard colorable one by adding a half twist to each loop. See Figure \ref{12}.
\end{example}

\begin{lemma}\label{le3}
Let $G$ be an embedded graph and let $D$ be the set of $d$-edge of $G$ arising from an all-crossing direction of $G_{m}$. Then $G-D$ and $G^{*}-D^{c}$ are both Eulerian.
\end{lemma}

\begin {figure}[!htbp]
\centering
\includegraphics[width=15cm]{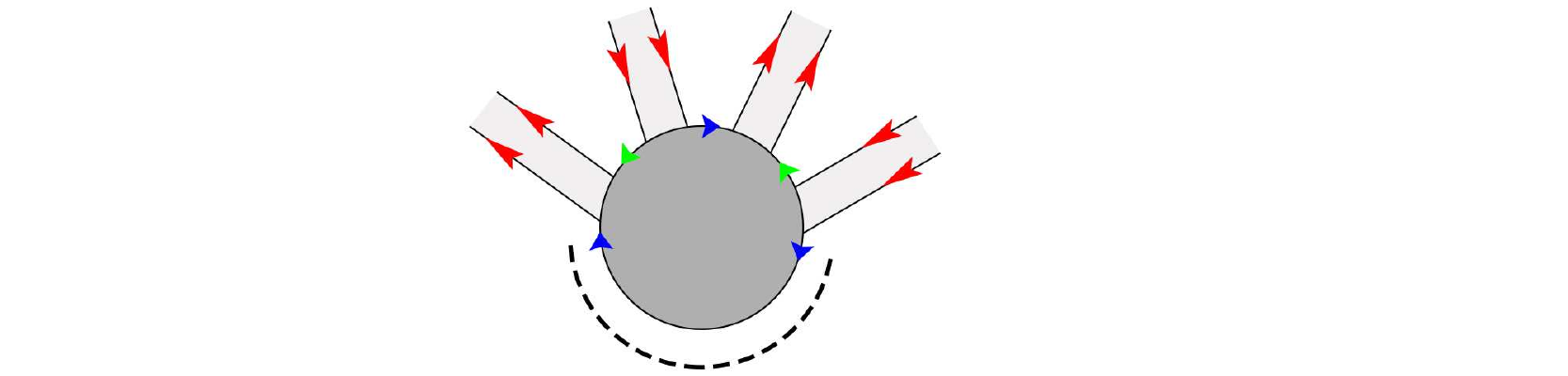}
\caption{Proof of Lemma \ref{le3}.}
\label{13}
\end {figure}

\noindent {\bf Proof.}
As shown in Figure \ref{8}, we obtain an orientation for each boundary component of $G-D$. The degree of each of the vertices of $G-D$ must be even as shown in Figure \ref{13}. Similarly, we can prove $G^{*}-D^{c}$ is also Eulerian.
~~

\vskip0.2cm

\begin{figure}[!htbp]
\begin{center}
\includegraphics[width=13cm]{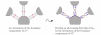}
\caption{Proof of Theorem \ref{th3}.}
\label{14}
\end{center}
\end{figure}

\noindent {\bf Proof of Theorem \ref{th3}.} $G^{\times}$ is orientable and assign it an anti-clockwise orientation. Then $G_m$ has an all-crossing direction as boundary components of $G^{\times}$. See Figure \ref{14}. Note that each edge of $G$ is a $d$-edge and apply Lemma \ref{le3}, the theorem is proved.

Let $G=(V(G),E(G))$ be a ribbon graph and $e\in E(G)$. We denote by $G-e$ the ribbon graph obtained from $G$ by deleting the edge $e$. We denote by $G/e=G^{e}-e$ the ribbon graph obtained from $G$ by contracting the edge $e$.  Let $A, B\subseteq E(G)$ and $A\cap B=\emptyset$. We denote by $G-A/B$ the ribbon graph obtained from $G$ by deleting the edges in $A$ and contracting the edges in $B$. We say that a ribbon graph $H$ is a {\it ribbon graph minor} of a ribbon graph $G$ if $H$ is obtained from $G$ by a sequence
of edge deletions, vertex deletions, or edge contractions.
See Moffatt\cite{MO} for the details of these constructions.
If $H$ is a ribbon graph minor of $G$, and $A\subseteq E(G)$, then by $H^{A}$ we mean $H^{A\cap E(H)}$.

\begin{lemma}\label{le1}
If $G$ is any embedded graph, $A, B, C\subseteq E(G)$ and $B\cap C=\emptyset$, then $(G-B/C)^{A}=G^{A}-B'/C'$, which $B'=(B\cap A^{c})\cup (C\cap A),$
$C'=(C\cap A^{c})\cup (B\cap A)$.
\end{lemma}

\noindent {\bf Proof.}
It is straightforward to check that $|E(G)|\in \{0, 1\}$ or $A=\emptyset$, so assume that $|E(G)|\geq 2$ and $A\neq\emptyset$.
By the definition of partial dual, it is enough to prove the lemma for $A=\{e\}$.
Let $f\in E(G)$ with $f\neq e$.
Since edge deletion, contraction and partial duality change the ribbon graph locally at the edge involved,
it follows that $(G-f)^{e}=G^{e}-f$ and $(G/f)^{e}=G^{e}/f$. By relation $G/e=G^{e}-e$ between partial duality and contraction,
$(G-e)^{e}=G-e=(G^{e})^{e}-e=(G^{e})/e$. Evidently, $(G/e)^{e}=G/e=G^{e}-e$.
It follows that $(G-B/C)^{A}=G^{A}-((B\cap A^{c})\cup (C\cap A))/((C\cap A^{c})\cup (B\cap A))=G^{A}-B'/C'$.

\begin{lemma}\label{le2}
If $G$ is any embedded graph, $A\subseteq E(G)$, and $G^{A}$ is bipartite, then $(G-{A^{c}})^{\ast}$ and $(G^{\ast}-{A})^{\ast}$ are both bipartite.
\end{lemma}

\noindent {\bf Proof.}
$(G-{A^{c}})^{\ast}=G^{*}/{A^{c}}=(G^{*})^{A^{c}}-A^{c}=G^{A}-A^{c}$, where the first equality is by Lemma \ref{le1},
the second by the relation between partial duality and contraction, and the third by the basic properties of partial duals.
Similarly, $(G^{\ast}-{A})^{\ast}=(G^{\ast})^{*}/A=G/A=G^{A}-A$. Since $G^{A}$ is bipartite, then $G^{A}-{A^{c}}$ and $G^{A}-A$ are bipartite, completing the proof.

\vskip0.2cm

Theorem \ref{th1} follows from Lemma \ref{le2} directly by the geometrical duality between the biparticity and checkerboard colourability of embedded graphs.

\begin{remark}
For plane graphs, Huggett and Moffatt \cite{HI} obtained the same results as Lemmas \ref{le3} and \ref{le2}.
\end{remark}

\begin {figure}[!htbp]
\centering
\includegraphics[width=12cm]{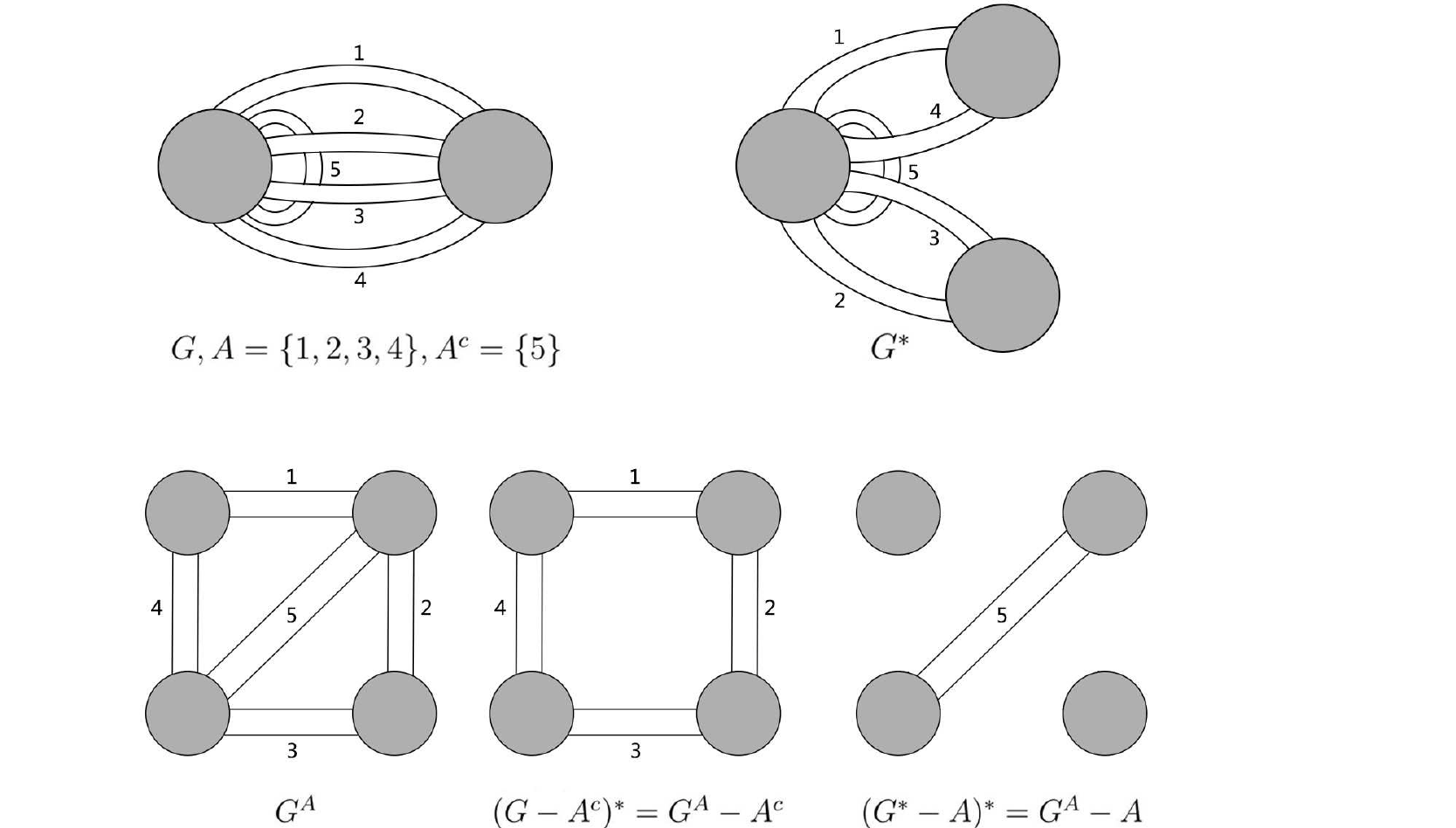}
\caption{A counterexample.}
\label{15}
\end {figure}

\begin{remark}
The converse of Lemma \ref{le2} is not true. A counterexample is given in Figure \ref{15}.
\end{remark}


\section*{Acknowledgements}
\noindent

This work is supported by NSFC (No. 11671336) and President's Funds of Xiamen University (No. 20720160011).

\section*{References}
\bibliographystyle{model1b-num-names}
\bibliography{<your-bib-database>}

\end{document}